\documentclass[12pt]{amsart}
\usepackage{a4wide}
\usepackage{amssymb}
\usepackage{amsthm}
\usepackage{amsmath}
\usepackage{amscd}
\usepackage[mathscr]{eucal}
\usepackage[all]{xy}

\usepackage{verbatim}

\numberwithin{equation}{section}

\theoremstyle{plain}
\newtheorem{theorem}{Theorem}[section]
\newtheorem{corollary}[theorem]{Corollary}

\newtheorem{proposition}[theorem]{Proposition}

\theoremstyle{definition}

\newtheorem{remark}[theorem]{Remark}

\theoremstyle{remark}

\newcommand{\R}{\mathbb{R}}
\newcommand{\Q}{\mathbb{Q}}
\newcommand{\Z}{\mathbb{Z}}

\newcommand{\C}{\mathbb{C}}
\newcommand{\h}{\mathbb{H}}
\renewcommand{\H}{\mathbb{H}}

\newcommand{\G}{\Gamma}
\newcommand{\g}{\gamma}

\newcommand{\back}{\backslash}

\newcommand{\zxz}[4]{\begin{pmatrix} #1 & #2 \\ #3 & #4 \end{pmatrix}}

\newcommand{\kzxz}[4]{\left(\begin{smallmatrix} #1 & #2 \\ #3 & #4\end{smallmatrix}\right) }
\newcommand{\kabcd}{\kzxz{a}{b}{c}{d}}

\newcommand{\calE}{\mathcal{E}}
\newcommand{\calF}{\mathcal{F}}

\newcommand{\calQ}{\mathcal{Q}}

\newcommand{\frake}{\mathfrak e}

\newcommand{\eps}{\varepsilon}
\newcommand{\bs}{\backslash}

\newcommand{\tr}{\operatorname{tr}}

\newcommand{\Span}{\operatorname{span}}

\newcommand{\Sl}{\operatorname{SL}}

\newcommand{\SL}{\operatorname{SL}}

\newcommand{\Mat}{\operatorname{Mat}}
\newcommand{\Spin}{\operatorname{Spin}}
\newcommand{\PSL}{\operatorname{PSL}}

\newcommand{\Orth}{\operatorname{O}}

\newcommand{\SO}{\operatorname{SO}}

\newcommand{\Iso}{\operatorname{Iso}}

\begin{document}

\title[Traces of CM values of modular functions]{Traces of CM values of modular functions and related topics}

\author[Jan H.~Bruinier]{Jan Hendrik Bruinier}

\address{Mathematisches Institut, Universit\"at zu K\"oln, Weyertal 86--90, D-50931 K\"oln, Germany}
\email{bruinier@math.uni-koeln.de }


\date{\today}

\begin{abstract}
The purpose of this note is to report on recent joint work with J.~Funke, P.~Jenkins, and K.~Ono on the traces of
CM values of modular functions and some applications \cite{BF}, \cite{BJO}.
\end{abstract}

\maketitle

\section{Introduction}

The classical  $j$-function on the complex upper half plane $\H$ is defined by
$$
j(\tau)=\frac{E_4(\tau)^3}{\eta(\tau)^{24}}= q^{-1}+744+196884q+21493760q^2+\dots.
$$
Here $\eta=q^{1/24}\prod_{n=1}^\infty (1-q^n)$ is the Dedekind eta function, $E_4=1+240\sum_{n=1}^{\infty}
\sum_{m\mid n}m^3q^n$ is the normalized Eisenstein series of weight $4$ for the group $\Gamma(1)=\PSL_2(\Z)$, and
$q=e(\tau)=e^{2\pi i \tau}$ for $\tau\in\H$. The $j$-function  is a Hauptmodul for the group $\Gamma(1)$, i.e., it
generates the field of all meromorphic modular functions for this group.

The values of $j(\tau)$ at CM points are known as {\it singular moduli}. They are algebraic integers generating
Hilbert class fields of imaginary quadratic fields. In this note we consider the traces of singular moduli and
more generally the traces of CM values of modular functions on modular curves of arbitrary genus.

Let $D$ be a positive integer and write $\calQ_{D}$ for the set of positive definite integral binary quadratic
forms $[a,b,c]$ of discriminant $-D=b^2-4ac$. The group $\Gamma(1)$ acts on $\calQ_{D}$. If $Q=[a,b,c]\in
\calQ_{D}$ we write $\Gamma(1)_Q$ for the stabilizer of $Q$ in $\Gamma(1)$ and
$\alpha_Q=\frac{-b+i\sqrt{D}}{2a}$ for the corresponding CM point in $\H$. By the theory of complex multiplication,
the values of $j$ at such points $\alpha_Q$ are algebraic integers whose degree is equal to the class number of
$K=\Q(\sqrt{-D})$. Moreover, $K(j(\alpha_Q))$ is the Hilbert class field of $K$. In \cite{GZ}, Gross and Zagier
derived a closed formula for the norm to $\Z$ of $j(\alpha_Q)$ as a special case of their work on the Gross-Zagier
formula. In a later paper \cite{Za2}, Zagier studies the trace of $j(\alpha_Q)$. We briefly recall his result.

To this end it is  convenient to consider the normalized Hauptmodul $J(\tau)=j(\tau)-744$ for $\Gamma(1)$ instead
of $j(\tau)$ itself. The modular trace of $J$ of index $D$ is defined as
\begin{equation}\label{modtrace}
{\bf t}_J(D) = \sum_{Q\in \calQ_{D}/\Gamma(1)} \frac{1}{|\Gamma(1)_Q|} J(\alpha_Q).
\end{equation}
Zagier discovered that the generating series
\begin{equation}\label{Zagier1}
-q^{-1} +2 + \sum_{D=1}^{\infty} {\bf t}_J(D) q^D= -q^{-1} +2 -248 q^3+492 q^4-4119q^7+7256q^8+\dots
\end{equation}
is a meromorphic modular form of weight $3/2$ for the Hecke group $\G_0(4)$ whose poles are supported at the
cusps. More precisely, it is equal to the weight $3/2$ form
\begin{equation}\label{Zagier2}
g(\tau)=\frac{\eta(\tau)^2 E_4(4\tau)}{\eta(2\tau)\eta(4\tau)^6}.
\end{equation}
Zagier gives two different proofs of this result. The first uses certain recursion relations for the ${\bf
t}_J(D)$, the second uses Borcherds products on $\Sl_2(\Z)$ and an application of Serre duality. Both proofs rely
on the fact that (the compactification of) $\G(1) \back \h$ has genus zero. In  \cite{Kim1,Kim2}, Kim extends
Zagier's results to other modular curves of genus zero using similar methods.

The above connection between the weight $3/2$ form $g$ for
$\Gamma_0(4)$ and the weight $0$ form $J$ for $\Gamma(1)$ reminds us
of (a special case of) the Shimura lift which is a linear map from
holomorphic modular forms of weight $k+1/2$ for $\Gamma_0(4)$ in the
Kohnen plus space to holomorphic modular forms of weight $2k$ for
$\Gamma(1)$. Moreover, it reminds of the theta lift from weight
$0$ Maass wave forms to weight $1/2$ Maass forms first considered by
Maass \cite{Maass} and later reconsidered by Duke \cite{Du1} and Katok
and Sarnak \cite{KS}.  However, there are two obvious differences:
First, in our case the half integral weight form has weight $3/2$
rather than $1/2$; and second, neither $J$ nor $g$ is holomorphic at
the cusps. The first difference should be not so serious, since there
is often a duality between weight $k$ and weight $2-k$ forms on
modular curves as a consequence of Serre duality. If we ignore the
second difference for a moment, in view of the work of Shintani
\cite{Sh} and Niwa \cite{Ni} realizing the Shimura lift as a theta
lift, it is natural to ask, whether Zagier's result can also be
interpreted in the light of the theta correspondence?

In other words, one might ask if there is a suitable theta function
$\theta(\tau, z,\varphi)$ which transforms like a modular form of
weight $3/2$ in $\tau$ and is invariant under $\Gamma(1)$ in $z$ such
that $g(\tau)$ is equal to the theta integral
\begin{align}
I(\tau,J)= \int_{\Gamma(1)\bs\H} J(z)\theta(\tau,z,\varphi) \,\frac{dx\,dy}{y^2}.
\end{align}
Clearly one has to be very careful with the convergence of the
integral because of the pole of $J$ at the cusp. It is shown in
\cite{BF} that it is possible to obtain such a description by
considering the theta kernel corresponding to a particular Schwartz
function $\varphi$ constructed by Kudla and Millson \cite{KM1}. This
generalizes \cite{Fu} where the lifting $I(\tau,1)$ of the constant
function $1$ was studied. A very nice feature of the theta kernel is
its very rapid decay at the cusps which leads to absolute convergence
of the integral.

The theta lift description of the correspondence between $J$ and $g$
can now be used to generalize Zagier's result to modular functions
(with poles supported at the cusps) on modular curves of arbitrary
genus. We will discuss this in section \ref{sect:2}. Moreover, one can
study the lifting for other automorphic functions. It turns out that
already the lifting of the non-holomorphic Eisenstein series
$E_0(z,s)$ of weight $0$ for $\Gamma(1)$ provides interesting
geometric and arithmetic insights. We will discuss this in section
\ref{sect:3}.

One can use the modularity of the generating series for the traces of
CM values of a modular function to obtain exact formulas. In section
\ref{sect:4} we briefly report on results of \cite{BJO} and \cite{Du2}
giving exact formulas for ${\bf t}_J(D)$ analogous to the
Hardy-Rademacher-Ramanujan formula for the partition function $p(n)$.
Moreover, we discuss the asymptotic behavior of ${\bf t}_J(D)$ as $D$
goes to infinity.


\section{The theta lift}
\label{sect:2}

Here we describe the theta lift discussed in the introduction. As in
the case of Shintani and Niwa it uses the dual pair
$\widetilde{\Sl}_2(\R)$, $\SO(1,2)$, combined with an exceptional
isomorphism relating $\SO(1,2)$ and $\Sl_2(\R)$.

\subsection{Setting}

Let $(V,q)$ be the quadratic space over $\Q$ of signature $(1,2)$
given by the trace zero $2\times 2 $ matrices,
\begin{equation} \label{iso}
V  :=\left\{ X = \begin{pmatrix} x_1 & x_2 \\ x_3 & -x_1
 \end{pmatrix} \in \Mat_2(\Q) \right\},
\end{equation}
with the quadratic form $q(X) = \det(X)$. The corresponding bilinear
form is $(X,Y) = -\tr(XY)$. (Note: for simplicity we assume that the
discriminant $d$ of the quadratic space is $1$. The more general case
is considered in \cite{BF}.) We let $G = \Spin(V)\simeq \SL_2$, viewed
as an algebraic group over $\Q$, and write $\bar{G}\simeq \PSL_2$ for
the image in $\SO(V)$. We realize the associated
symmetric space $D$ as 
the Grassmannian of lines in $V(\R)$ on which the quadratic form $q$ is
positive definite:
\[
D \simeq \{ z \subset V(\R) ; \;\text{$\dim z =1$ and $q|_z > 0$} \}.
\]

The group $\Sl_2(\Q)$ acts on $V$ by conjugation,
\[
g.X := gXg^{-1}
\]
for $X\in V$ and $g\in \Sl_2(\Q)$. This gives rise to an isomorphism 
$G\simeq\SL_2$.

Moreover, $D$ can be identified with the complex upper half plane $\H$
as follows: We pick as a base point $z_0\in {D}$ the line spanned by
$\left( \begin{smallmatrix} 0 & 1 \\ -1 & 0 \end{smallmatrix}\right)$,
and note that ${K} = \SO(2)$ is its stabilizer in ${G}(\R)$. For $z
\in \h$, we chose $g_z \in {G}(\R)$ such that $g_zi = z$.
We obtain the isomorphism $\h  \to
{D}$,
\begin{equation}
z \longmapsto g_z z_0 = \Span\left(g_z . \left( \begin{smallmatrix} 0 & 1 \\ -1 & 0
\end{smallmatrix}\right)\right).
\end{equation}
So for $z = x + iy \in \h$, the associated positive line is generated by
\begin{equation}
X(z) :=  g_z.\zxz{0}{1}{-1}{0} = \frac1{y}
\begin{pmatrix} -\frac12(z+ \bar{z}) & z\bar{z} \\ -1 & \frac12(z+ \bar{z})
\end{pmatrix}.
\end{equation}
In particular, $q(X(z)) = 1$ and $g.X(z) = X(gz)$ for $g \in G(\R)$.

Let $L \subset V$ be an even lattice of full rank and write $L^\#$ for the dual lattice of $L$. Let $\G$ be a
congruence subgroup of $\Spin(L)$ which takes $L$ to itself and acts trivially on the discriminant group
$L^{\#}/L$. We write $M = \G \backslash D$ for the associated modular curve.

We now define CM points in this setting. For $X \in V(\Q)$ of positive norm we put
\begin{equation}
D_X = \Span(X) \in D.
\end{equation}
It is easily seen that the corresponding point in $\H$ satisfies a quadratic equation over $\Q$. The stabilizer
$G_X$ of $X$ in $G(\R)$ is isomorphic to $\SO(2)$ and $\G_X = G_X \cap \G$ 
is finite.
For $m \in \Q_{>0}$ and a congruence condition $h \in L^{\#}$,  the group $\G$ acts on
\[
L_{h,m} = \{X \in L +h ;\; q(X) =m\}
\]
with finitely many orbits. We define the \emph{Heegner divisor} of discriminant $m$  on $M$ by
\begin{equation}
Z(h,m) = \sum_{X \in \G \back L_{h,m} }\frac{1}{|\bar{\G}_X|} D_X.
\end{equation}

\subsection{The Kudla Millson theta series}

In \cite{KM1}, Kudla and Millson constructed (in greater generality) a Schwartz function $\varphi$ on $V(\R)$ with
values in $\Omega^{1,1}(D)$, the closed differential forms on $D$ of Hodge type $(1,1)$. In our particular case
it is given by
\begin{equation}
\varphi(X,z) = \biggl( (X,X(z))^2  - \frac1{2\pi} \biggr) \, e^{-\pi (X,X(z))^2 +\pi (X,X)} \, \mu,
\end{equation}
where $X\in V(\R)$, $z=x+iy\in \H$, and $ \mu = \tfrac{dx \wedge
dy}{y^2} = \tfrac{i}2 \tfrac{dz \wedge d\bar{z}}{y^2}$. Notice that
$\varphi (g.X,g z) = \varphi (X,z)$ for $g \in G(\R)$. We
put
\begin{align}
\varphi^0(X,z) = e^{\pi(X,X)} \varphi (X,z) =\left((X,X(z))^2  - \frac1{2\pi} \right) \, e^{-\pi (X,X(z))^2 +2\pi
(X,X)} \, \mu.
\end{align}
The geometric significance of this Schwartz function lies in the fact
that for $q(X) >0$, the $2$-form $\varphi^0(X,z)$ is a Poincar\'e dual
form for the CM point $D_X$, while it is exact for $ q(X)
<0$. 

As usual, from the Schwartz function $\varphi$ one can construct a
theta series as follows. We let $\omega$ be the Weil representation of
$\widetilde{\Sl}_2(\R)$ on the Schwartz space associated to the
additive character $t \mapsto e^{2\pi i t}$. For $\tau = u+ iv \in
\h$, we put $g'_{\tau} = \left(
\begin{smallmatrix}1&u\\0&1\end{smallmatrix} \right) \left(
\begin{smallmatrix}v^{1/2}&0\\0&v^{-1/2}\end{smallmatrix} \right)$, 
so that $g'_{\tau} i =\tau$, and define
\begin{align*}
\varphi(X,\tau,z) := v^{-3/4} \omega( g'_{\tau}) {\varphi}(X,z) =
e^{2\pi i q(X) \tau} \varphi^0(\sqrt{v}X,z).
\end{align*}
Then,  for $h \in L^{\#}/L$, the theta kernel
\begin{equation}
\theta_{h}(\tau,z,\varphi) = \sum_{X \in h + L} \varphi(X,\tau,z)
\end{equation}
has a nice transformation behavior in both variables, $\tau$ and $z$
(see \cite{KM3}, \cite{Fu}). It is a $\Gamma$-invariant differential form in
$z$, and transforms as a \emph{non-holomorphic} modular form of weight
$3/2$ for the congruence subgroup $\G(N)$ of $\SL_2(\Z)$, where $N$ is
the level of the lattice $L$. 
To lighten
the notation, we will frequently drop the argument $\varphi$.

A rather surprising and important feature of this theta series is its rapid decay at the boundary:

\begin{proposition}[\cite{Fu}, Proposition 4.1]\label{thetagrowth}
Write $z=x+iy$ with $x,y\in \R$ and let $\sigma\in \Gamma(1)$. There
is a constant $C>0$ such that
\[
\theta_{h}(\tau,\sigma z) = O(e^{-Cy^2}), \qquad \text{$y \to \infty$},
\]
uniformly in $x$.
\end{proposition}

\subsection{The theta integral}

We denote by $ M^!_0(\G)$ the space of (scalar valued) weakly
holomorphic modular forms of weight $0$ with respect to $\G$. It
consists of those meromorphic modular functions for $\G$ which are
holomorphic outside the cusps. So, for instance, $M^!_0(\Sl_2(\Z))=\C[J]$.

If $h\in L^\#/L$ and $f\in M^!_0(\G)$, we define a theta integral by
\begin{equation}\label{comp-int}
I_h(\tau,f) = \int_M f(z) \theta_{h}(\tau,z).
\end{equation}
Proposition~\ref{thetagrowth} implies the convergence of the integral,
since the decay of $\theta_{h}(\tau,z)$ offsets the exponential growth
of $f$ at the cusps. Consequently, $I(\tau,f)$ defines a (in general
non-holomorphic) modular form of weight $3/2$.

Now the main task is to determine the Fourier expansion of
$I_h(\tau,f)$. The computation of the Fourier coefficients with
positive index is quite straightforward. Using the properties of the
Kudla-Millson Schwartz function $\varphi$, it can be shown that they
are given by traces of CM values of $f$. However, the constant term of
$I_h(\tau,f)$ and the negative coefficients are more involved. Here
convergence becomes a subtle issue. These calculations are the
technical heart of \cite{BF}. Eventually they lead to the following
theorem (cf.~\cite{BF}, Theorem~4.5, Proposition~4.7, Corollary~4.8).

\begin{theorem}\label{MAIN}
Let $f \in M_0^!(\G)$ and assume that the constant coefficients of $f$
at all cusps of $M$ vanish. Then $I_h(\tau,f)$ is a weakly holomorphic
modular form of weight $3/2$ for $\Gamma(N)$. The Fourier expansion of
$I_h(\tau,f)$ is given by
\begin{align*}
I_h(\tau,f) &=   \sum_{\substack{m\in \Z+q(h)\\ m\gg-\infty}} {\bf{t}}_f(h,m) q^m,
\end{align*}
where ${\bf{t}}_f(h,m)$ is the modular trace function,
\begin{align*}
{\bf{t}}_f(h,m) =
\begin{cases}
\displaystyle
\sum_{X \in \G \back L_{h,m}} \frac1{|\bar{\G}_X|} f(D_X),&\text{if $m>0$,}\\[3ex]
\displaystyle
-  \frac{\delta_{h,0}}{2\pi} \int^{reg}_M f(z) \,\frac{dx\,dy}{y^2},&\text{if $m=0$,}\\[3ex]
\displaystyle \text{\parbox{8cm}{explicit formula in terms of geodesic cycles connecting two cusps,}} &\text{if $m<0$.}
\end{cases}
\end{align*}
Here $\delta_{h,0}$ denotes the Kronecker delta, and for the precise definition of ${\bf{t}}_f(h,m)$ for $m<0$ we
refer to \cite{BF} Definition~4.4. If the constant coefficients of $f$ do not vanish, then in addition, the Fourier
expansion contains certain non-holomorphic terms which are supported on Fourier coefficients with non-positive
index.
\end{theorem}

The regularized integral occurring in the constant term is defined as
\[
\lim_{\eps\to 0}\int_{M(\eps)} f(z) \,\frac{dx\,dy}{y^2},
\]
where $M(\eps)$ denotes the manifold with boundary obtained by removing an $\eps$-disc around each cusp from $M$.
It can be viewed as a regularized average value of $f$, and it can be explicitly computed by \cite{BF} Remark~4.9.

\begin{remark}\label{remark2.3}
Let $\widetilde{\Sl}_2(\R)$ be the metaplectic two-fold cover of $\SL_2(\R)$ realized by the two choices of
holomorphic square roots of $\tau \mapsto  c\tau + d$ for $g = \left(
\begin{smallmatrix} a&b \\ c&d \end{smallmatrix} \right) \in
\SL_2(\R)$. Recall that there is a unitary representation $\rho_L$ of the inverse image $\G'$ of $\Sl_2(\Z)$ under the
covering map, acting on the group algebra $\C[L^{\#}/L]$ (see \cite{Bo1}, \cite{Br2}). We denote the standard
basis elements of $\C[L^{\#}/L]$ by $\frake_{h}$, where $h \in {L^{\#}/L}$. For the generators $S = \left( \left(
\begin{smallmatrix} 0&-1\\1&0
\end{smallmatrix} \right), \sqrt{\tau} \right)$
and $T=\left( \left( \begin{smallmatrix} 1 &1 \\ 0&1 \end{smallmatrix} \right),1\right)$ of $\G'$ the action of
$\rho_L$ is given by
\begin{align*}
\rho_L(T) \frake_{h} &= e((h,h)/2) \frake_{h},\\
\rho_L(S) \frake_{h} &= \frac{\sqrt{i}}{\sqrt{|L^{\#}/L|}} \sum_{h'\in L^{\#}/L} e(- (h,h')) \frake_{h'}.
\end{align*}
If one considers instead of the individual integral $I_h(\tau,f)$ the vector
\begin{align}
I(\tau,f)=\sum_{h\in L^\#/L} I_h(\tau,f) \frake_h,
\end{align}
one obtains a vector valued modular form of weight $3/2$ for the group $\G'$ and the representation $\rho_L$. If
the discriminant group $L^\#/L$ of $L$ is cyclic then such vector valued modular forms can also be interpreted as
weak Jacobi forms in the sense of \cite{EZ}.
\end{remark}

We end this section with an example illustrating the theorem. Let $p$ be a prime (or $p=1$), and let $L$ be the
lattice
\[
L=\left\{\zxz{b}{2c}{2ap}{-b};\quad a,b,c\in \Z\right\}\subset V.
\]
Then $L$ has level $4p$ and is stabilized by $\G_0^*(p)$, the extension of the Hecke group
$\Gamma_0(p)\subset\Gamma(1)$ with the Fricke involution $W_p=\kzxz{0}{-1}{p}{0}$. We take $\G=\G_0^*(p)$ so that
$M$ is the modular curve $\G_0^*(p)\bs\H$.

For a positive integer $D$, we consider the subset $\calQ_{D,p}$ of quadratic forms $[a,b,c]\in \calQ_{D}$ such
that $a\equiv 0\pmod{p}$. The group $\G_0^{\ast}(p)$ acts on $\calQ_{D,p}$ with finitely many orbits. It turns out
that the Heegner divisor $Z(0,D)$ on $M$ is equal to
\begin{align}
 \sum_{Q\in \calQ_{D,p}/\Gamma_0^*(p)}\frac{1}{|\Gamma_0^*(p)_Q|}\alpha_Q,
\end{align}
where $\Gamma^{\ast}_0(p)_Q$ is the stabilizer of $Q$ in $\Gamma^{\ast}_0(p)$. Consequently, if  $f$ is a weakly
holomorphic modular function (of weight $0$) for $\Gamma^{\ast}_0(p)$, then the modular trace of $f$ over the
Heegner divisor $Z(0,D)$ can also be written as
\begin{equation}\label{modtrace2}
{\bf t}^{\ast}_f(D) = \sum_{Q\in \calQ_{D,p}/\Gamma^{\ast}_0(p)} \frac{1}{|\Gamma^{\ast}_0(p)_Q|} f(\alpha_Q).
\end{equation}

\begin{theorem} \label{thm:intro}
Let $f=\sum_{n\gg-\infty}a(n)q^n\in M^!_0(\Gamma_0^*(p))$ and assume that the constant coefficient $a(0)$ vanishes.
Then
\begin{align*}
\frac{1}{2}I_0(\tau,f)&=\sum_{\substack{D>0}}  {\bf t}^{\ast}_f(D) q^D +
\sum_{n\geq 0}\big(\sigma_1(n)  + p \sigma_1(n/p) \big) a(-n) \\
& \phantom{=}{}-\sum_{m>0}\sum_{n>0}m a(-mn)q^{-m^2}
\end{align*}
is a weakly holomorphic  modular form of weight $3/2$ for the group $\G_0(4p)$ satisfying the Kohnen plus space
condition. Here $\sigma_1(0)=-\frac{1}{24}$ and $\sigma_1(n)=\sum_{t\mid n} t$ for $n\in \Z_{\geq0}$ and
$\sigma_1(x) =0$ for $x \notin \Z_{\geq0}$.


\end{theorem}

For $p=1$, and  $f=J$, we recover Zagier's result \eqref{Zagier1}.

\section{Extensions}
\label{sect:3}

It is natural to consider the theta lift of the previous section for other automorphic functions. Already the
lifting of the real analytic Eisenstein series of weight $0$ for $\Sl_2(\Z)$ turns out to be quite interesting. We
normalize this Eisenstein series as follows:
\[
\calE_0(z,s) = \frac12 \zeta^{\ast}(2s+1) \sum_{\g \in \G_{\infty} \back \Sl_2(\Z)} \left(\Im(\g
z)\right)^{s+\frac12}.
\]
Here $\G_{\infty} = \left\{\kzxz{1}{n}{0}{1};\; n\in \Z\right\}$ and $\zeta^{\ast}(s) = \pi^{-s/2} \G(s/2)
\zeta(s)$ is the completed Riemann zeta function. Recall that $\calE_0(z,s)$ converges for $\Re(s) >1/2$ and has a
meromorphic continuation to $\C$ with a simple pole at $s = 1/2$ with residue $1/2$. It satisfies the functional
equation $\calE_0(z,-s) = \calE_0(z,s)$.

We consider the lattice
\[
L=\left\{\zxz{b}{c}{a}{-b};\quad a,b,c\in \Z\right\}\subset V.
\]
We have $L^\#/L\cong \Z/2\Z$, the level of $L$ is $4$, and $\Gamma=\Sl_2(\Z)$ takes $L$ to itself and acts
trivially on $L^\#/L$. We let $\mathfrak{e}_0, \mathfrak{e}_1$ be the standard basis of  $\C[L^\#/L]$
corresponding to the cosets $\left( \begin{smallmatrix} b & 0 \\ 0& -b \end{smallmatrix} \right)$ with $b = 0$ and
$b = 1/2$, respectively.

We define a vector valued Eisenstein series $\calE_{3/2,L}(\tau,s)$ of weight $3/2$ for the Weil representation
$\rho_{L}$ (see Remark \ref{remark2.3}) by
\[
\calE_{3/2,L}(\tau,s) = -\frac1{4\pi} (s+\frac12) \zeta^{\ast}(2s+1) \sum_{\g' \in \G'_{\infty} \back \G'} \left(
\Im(\tau)^{\frac12(s-\frac12)} \mathfrak{e}_0\right)|_{3/2,L} \,\g'.
\]
Here the Petersson slash operator is defined on functions $f: \h \to
\C[L^{\#}/L]$ by
\[
(f |_{3/2,L}\, \g') (\tau) =  \phi(\tau)^{-3} \rho^{-1}_{L}(\g') f(\g \tau)
\]
for $\g'=(\g,\phi)\in \G'$. Moreover, $\G_{\infty}'$ is the inverse
image of $\G_{\infty}$ inside $\G'$. The argument of \cite{EZ} \S5
Theorem 5.4 implies that the scalar valued function
\begin{equation}\label{ZEisen}
\mathcal{F}(\tau,s) =  \left(\calE_{3/2,L}(4\tau,s) \right)_{0} + \left( \calE_{3/2,L}(4\tau,s) \right)_{1}
\end{equation}
is a non-holomorphic modular form of weight $3/2$ for $\Gamma_0(4)$
satisfying the Kohnen plus space condition. Up to a constant factor
depending only on $s$ it is equal to Zagier's Eisenstein series as in
\cite{HZ}, \cite{Za1}.  Note that our $\mathcal{F}(\tau,s)$ is equal
to the Eisenstein series $\mathcal{E}(\tau,s)$ of \cite{Yang} formula
(3.9).

By applying a partial Fourier transform to the theta kernel $\theta_h(\tau,z)$ and unfolding against the resulting
Poincar\'e type series one obtains the following theorem.

\begin{theorem}\label{Eisensteinlift}
The theta integral of $\calE_0(z,s)$ is given by
\begin{equation}\label{th71}
I(\tau,\calE_0(z,s)) =  \zeta^{\ast}(s+1/2) \calE_{3/2,L}(\tau,s).
\end{equation}
\end{theorem}

As a corollary one obtains another proof of the functional equation $\calE_{3/2,L}(\tau,-s)=
\calE_{3/2,L}(\tau,s)$.
Taking residues at $s=1/2$ on both sides of \eqref{th71} we obtain a different proof of Theorem~1.1 of \cite{Fu}:

\begin{corollary}\label{ZagierSW}
The theta integral of the constant function $1$ is given by
\[
I(\tau,1) =  2 \calE_{3/2,L}\left(\tau,1/2\right).
\]
\end{corollary}

On the other hand, the computation of the Fourier expansion of $I(\tau,1)$ 
in Theorem \ref{MAIN} shows that
\begin{equation}\label{Eisen-value}
\frac{1}{2}(I_0(\tau,1)+I_1(\tau,1))
=\sum_{D =0}^{\infty}  {\bf t}_1(D)   q^D + \frac{1}{16 \pi\sqrt{v}} \sum_{N=- \infty}^{\infty} \beta(4\pi N^2 v) q^{-N^2}.
\end{equation}
Here the modular trace ${\bf t}_1(D)$ over the Heegner divisor of discriminant $D>0$ is simply the Hurwitz-Kronecker class number $H(D)= \sum_{Q\in \calQ_{D}/\Gamma(1)} \frac{1}{|\Gamma(1)_Q|}$, and 
 ${\bf t}_1(0) = -\tfrac1{12}$.
Moreover, $\beta(s) = \int_1^{\infty} t^{-3/2} e^{-st} dt$.
By Corollary \ref{ZagierSW}, we find that \eqref{Eisen-value} is the Fourier expansion of $\mathcal{F}(\tau,1/2)$. Of course, the Fourier expansion of this Eisenstein series can also be computed directly (see \cite{Za1}). However, Theorem \ref{MAIN} provides a conceptual explanation of the geometric interpretation of the positive coefficients.
Notice that the theta integral $I(\tau,1)$ is a \emph{non-holomorphic} modular form. The non-holomorphic part supported on the non-positive Fourier coefficients is the prototype of the non-holomorphic contributions that occur in Theorem \ref{MAIN} if the constant coefficients of the input form $f$ do not all vanish. 

Now we consider the constant terms in the Laurent expansions on both sides of \eqref{th71}. The constant term of the Eisenstein series $\calE_0(z,s)$ is given by the Kronecker limit formula. We have
\begin{equation}\label{Kroni}
- \frac{1}{12}\log\left( |\Delta(z)y^6| \right) = \lim_{s \to 1/2} \left( \calE_0(z,s) - \zeta^{\ast}(2s-1)\right),
\end{equation}
where $\Delta(z)=\eta(z)^{24}$ is the classical Delta function.
On the right hand side the constant term will involve the derivative 
$\calE_{3/2,L}'(\tau,1/2)$. Working out the details we obtain:

\begin{theorem}
Putting
$\|\Delta(z)\| = e^{-3(\gamma+\log(4\pi))} |\Delta(z) (4\pi y)^6|$ 
as in \cite{Yang},
we have
\[
- \frac{1}{12} I\left(\tau,\log \|\Delta(z)\| \right) = \calE'_{3/2,L}\left(\tau,1/2\right).
\]
\end{theorem}

Again, using the properties of the Kudla-Millson Schwartz function
$\varphi$ and the corresponding Green function $\xi$ constructed by
Kudla \cite{KAnn}, one can obtain a geometric interpretation of the
Fourier coefficients of the theta integral. It turns out that the
$D$-th coefficient of $-\frac1{12} I(\tau, \log\|\Delta(z)\|)$
is equal to an arithmetic intersection pairing $4 \langle
\widehat{\mathcal{Z}}(D,v), \widehat{\omega}\rangle$ in the sense of 
\cite{Bost,Kuehn,Soule}. Here
$\widehat{\omega}$ is the normalized metrized Hodge bundle on the
moduli stack over $\Z$ of elliptic curves.
Moreover, $\widehat{\mathcal{Z}}(D,v)$ is the arithmetic divisor 
given by the Heegner points of discriminant $D$ over $\Z$ 
and the corresponding
Kudla Green function at the archimedian place (see \cite{Yang}, \cite{BF} for
details).
In that way, one obtains a somewhat more direct proof of
the result of \cite{Yang}, stating that
$\frac{1}{4}\calF'(\tau,1/2)$ is the generating series for the
arithmetic degrees $\langle \widehat{\mathcal{Z}}(D,v),
\widehat{\omega}\rangle$. It will be interesting to extend this
argument to modular curves of arbitrary level.

\section{Exact formulas and asymptotics}
\label{sect:4}

In this section we come back to the modular curve $\Gamma(1)\bs\H$ 
and consider mainly the modular function $J$. 
We briefly describe how the modularity of the generating series 
for the traces of singular moduli  can be used to obtain 
exact formulas for ${\bf t}_J(D)$. 
Moreover, we describe some asymptotic results.
For more general results in this direction we refer to  \cite{Du2}, \cite{BJO}.

We write $M^{!,+}_{3/2}$ for the space of weakly holomorphic modular 
forms of weight $3/2$ for the group $\Gamma_0(4)$ satisfying 
the Kohnen plus space condition.
Recall from Theorem \ref{thm:intro} that for 
$f\in M^!_0(\Gamma(1))=\C[J]$ with Fourier expansion 
$f=\sum_{n\gg-\infty}a(n)q^n$ and $a(0)=0$, the generating series
\[
\frac{1}{2}I_0(\tau,f)=\sum_{\substack{D>0}}  {\bf t}_f(D) q^D +2
\sum_{n\geq 0}\sigma_1(n)  a(-n) -\sum_{m>0}\sum_{n>0}m a(-mn)q^{-m^2}
\]
belongs to $M^{!,+}_{3/2}$.
In particular, for the case $f=J$ considered in the introduction, the weight $3/2$ form is explicitly given by  
\eqref{Zagier2}.

Recall that the generating series for the classical partition function $p(n)$, 
\[
\eta(\tau)^{-1} = q^{-1/24}\sum_{n=0}^\infty p(n) q^n,
\]
is a modular form of weight $-1/2$ for the group $\Gamma(1)$ with a
multiplier system which can be described in terms of Dedekind sums.
This fact was used by Hardy, Ramanujan, and Rademacher, to obtain a
closed formula for $p(n)$ as an infinite series by means of the circle
method (see \cite{Ap} chapter 5).  Hejhal pointed out  that one can
use non-holomorphic Poincar\'e series to give a 
somewhat more conceptual proof of this result (see \cite{He} pp.~654).
It is natural to apply similar arguments for the generating series 
$g(\tau)$ of the ${\bf t}_J(D)$. Theorem~1.2 of \cite{BJO} implies the following
result:

\begin{theorem}\label{mainthm2}
We have 
$$ 
{\bf t}_J(D)=-24H(D)+\sum_{\substack{c>0\\ c\equiv 0\; (4)}}
S(D,c)\sinh(4\pi\sqrt{D}/c), 
$$
where $S(D,c)$ is the exponential sum 
$$
S(D,c)=\sum_{x^2\equiv -D\;(c)}e(2x/c).
$$
\end{theorem}

We sketch the basic idea of the proof which is rather simple.
To avoid technical complications, instead of looking at
$\eta(\tau)^{-1}$ or $g(\tau)$ we first consider weakly holomorphic
modular forms of weight $k=4,6,8,10,14$ for the group $\Gamma(1)$ with
trivial multiplier system. We assume that the weight is greater than
$2$ to ensure absolute convergence of the Poincar\'e
series. Moreover, the upper bound on the weight ensures that there are
no cusp forms.

Let $m$ be a positive integer and let $f_m$ be the unique weakly
holomorphic form for $\Gamma(1)$ whose Fourier expansion has the form
\[
f_m=q^{-m} + O(q).
\]
It is easy to see that $f_m$ exists for every $m$. For $m=1$, one can
take $E_{k}(\tau)\cdot (j(\tau)+c)$, where $E_k$ is the normalized
Eisenstein series of weight $k$ for $\Gamma(1)$, and the constant $c$
is chosen such that the constant term vanishes. Now the $f_m$ can be
constructed inductively by multiplying $f_{m-1}$ with $j$ and
subtracting suitable multiples of $E_k, f_1,\dots,f_{m-1}$.  For
instance in weight $k=4$ we have
\begin{align*}
f_1 &= E_4(\tau)\cdot (j(\tau)-984)\\
&=q^{-1}+141444\cdot q+68234240\cdot q^2+6446476530\cdot 
q^3+275423256576\cdot q^4+\dots.
\end{align*}

In order to obtain a formula for the coefficients of $f_m$, we construct this modular form in a different way as a Poincar\'e series.
We consider
\begin{align}\label{fm}
F_m(\tau)= \sum_{\g \in \G_{\infty} \back \G(1)} q^{-m} |_k \g.
\end{align}
Here the Petersson slash operator is defined on functions $f: \h \to
\C$ by
\[
(f |_{k}\, \g) (\tau) =  (c\tau + d)^{-k} f(\g \tau)
\]
for $\g=\kabcd\in \G(1)$.
The series \eqref{fm} converges normally and therefore has the transformation behavior of a modular form of weight $k$. The trivial coset in the sum contributes the term $q^{-m}$. The rest of the sum decays as $\Im(\tau)\to \infty$.
Consequently, $F_m$ is a weakly holomorphic modular form with Fourier expansion $F_m=q^{-m}+O(q)$. Hence 
$F_m=f_m$.

Now the Fourier expansion of $F_m$ can be computed in the same way as for the usual holomorphic Poincare series.
If we write $F_m=\sum_{n\gg-\infty}a(n)q^n$, we find for $n>0$ that
\begin{align}
a(n)&=2\pi (-1)^{k/2} \left(\frac{n}{m}\right)^{\frac{k-1}{2}} \sum_{c=1}^\infty \frac{1}{c} K(m,n,c) \,I_{k-1}\!\left(\frac{4\pi}{c}\sqrt{mn}\right),
\end{align}
where $K(m,n,c)$ denotes the Kloosterman sum
\begin{align}
K(m,n,c)&=\sum_{d\;(c)^*} e\left(\frac{m \bar d + n d}{c}\right).
\end{align}
Here the sum runs through the primitive residues modulo $c$, and $\bar d$ denotes the multiplicative inverse of $d$ modulo $c$. Moreover $I_\nu$ is the usual Bessel function as in \cite{AbSt} \S9.

This is the easiest instance of a Hardy-Ramanujan-Rademacher type
formula for the coefficients of a weakly holomorphic modular form.  If
one tries to apply this argument to our generating series $g(\tau)$
several complications arise.  The function $g(\tau)$ is a modular form
only for the group $\G_0(4)$ which has three cusps. Only those linear
combinations of Poincar\'e series are relevant which belong to the
Kohnen plus space. This is a rather technical difficulty which can be
handled by looking at the Poincar\'e series at the cusp $\infty$ and
then applying the Kohnen projection operator to the plus space.
A more serious problem is that Poincar\'e series in weight
$3/2$ do not converge and have to be defined by analytic continuation
(``Hecke summation''). This can be done using the spectral theory of
the resolvent kernel \cite{Fa}, \cite{He}.

In that way, for every positive integer $m$, we obtain a Poincar\'e series $F_m^+(\tau)$, which transforms like a modular form of weight $3/2$ for $\G_0(4)$ and satisfies the plus space condition.
However, there is no reason for $F_m^+(\tau)$ to be holomorphic in $\tau$ as a function on $\H$, and it turns out that   $F_m^+(\tau)$ is in fact often {\em non-holomorphic}.

The good thing is that the non-holomorphic part can be computed explicitly. 
If $m$ is a square, one finds that (up to a constant multiple) 
it is the same as the non-holomorphic part of Zagier's Eisenstein series $\calF(\tau,1/2)$, see \eqref{Eisen-value}.
More precisely, 
\begin{align}\label{holsplit}
F_m^+(\tau)+24 \calF(\tau,1/2)=q^{-m} + O(1) \in M^{!,+}_{3/2}
\end{align}
is a weakly holomorphic modular form.
If $m$ is not a square, one finds that $F_m^+(\tau)\in M^{!,+}_{3/2}$ is already weakly holomorphic.

Now, from \eqref{holsplit} we easily deduce that 
\[
g(\tau)=-F_1^+(\tau)-24 \calF(\tau,1/2).
\]
The Fourier expansion of $F_1^+(\tau)$ can be computed and the expansion of $\calF(\tau,1/2)$ is given in \eqref{Eisen-value}. This leads to the formula of Theorem 
\ref{mainthm2}. The sum over $c>0$ comes from the $D$-th coefficient of $F_1^+$,
using the identity $I_{1/2}(z)=\left (\frac{2}{\pi z}\right )^{1/2} \sinh(z)$.
The class number comes from the $D$-th coefficient of $\calF(\tau,1/2)$.

It should be pointed out that the formula of Theorem \ref{mainthm2} is not very useful for numerical computations of  
${\bf t}_J(D)$, because of its slow rate of convergence.
However, it is possible to derive some asymptotic information from it.

Regarding the size of ${\bf t}_J(D)$, it is straightforward to see that for any $c>1/2$ we have
\[
{\bf t}_J(D) = (-1)^D e^{\pi\sqrt{D}} + O(e^{c\pi \sqrt{D}}),\quad D\to \infty.
\]
A much stronger result was conjectured (in a slightly different form) in \cite{BJO}
and recently proved by Duke \cite{Du2}. We close the present 
note by stating this result.

Recall that a positive definite integral binary quadratic form $Q\in\calQ_D$ is called {\em reduced}, if the corresponding point $\alpha_Q\in \H$ lies in the usual fundamental domain
\begin{equation}\label{fundomain}
\mathcal{F}=\left \{ -\frac{1}{2}\leq \Re(z) < \frac{1}{2}\
{\text {\rm and}}\  |z|>1 \right \} \cup \left \{-\frac{1}{2}\leq
\Re(z)\leq 0 \ \ {\text {\rm and}}\ \ |z|=1\right \}.
\end{equation}
We write $\calQ_D^{red}$ for the set of reduced quadratic forms in $\calQ_D$. So $\calQ_D^{red}$ is a set of representatives for $\calQ_D/\Gamma(1)$.

\begin{theorem}[Duke]\label{24thm}
As $-D$ ranges over negative fundamental discriminants, we have
\[ 
\lim_{D \rightarrow \infty} \frac{1}{H(D)}\Bigg( {\bf t}_J(D) -
\sum_{\substack{Q\in \calQ_D^{red}\\ \Im(\alpha_Q)>1}} e(-\alpha_Q)\Bigg) = -24.
\] 
\end{theorem}
So the finite sum over reduced quadratic forms on the left hand side describes the growth of  ${\bf t}_J(D)$ amazingly well.
The value $-24$ of the limit arises as the regularized average value 
\[
\frac{3}{\pi} \int_{\Gamma(1)\bs \H}^{reg} J(z) \,\frac{dx\,dy}{y^2}
\]
of $J$ which is also known as the Atkin functional on
$M_0^!(\Gamma(1))$.  Notice that Theorem \ref{24thm} holds in greater
generality for every $f\in M_0^!(\Gamma(1))$ with the appropriate
modifications, see \cite{Du2} Theorem~1.  
The proof relies on an  
application of the equidistribution of CM points \cite{Du1}.

Using Theorem \ref{mainthm2}, it is not difficult to show that
Theorem~\ref{24thm} is equivalent to the assertion that $$
\sum_{\substack{ c>\sqrt{D/3}\\ c\equiv 0\; (4) }}
S(D,c)\sinh\left (\frac{4\pi}{c} \sqrt{D}\right)
=o\left(H(D)\right). $$
In view of Siegel's theorem that $H(D)\gg_{\epsilon}
D^{\frac{1}{2}-\epsilon}$, Theorem \ref{24thm} would follow from a
bound for such sums of the form $\ll D^{\frac{1}{2}-\gamma}$, for some
$\gamma>0$.  Estimates of this type have been established for
such sums, but are difficult to establish.
They are implicit in Duke's proof of the
equidistribution of CM points and therefore in his proof of Theorem
\ref{24thm} as well.  
These bounds are intimately connected to the problem of
bounding coefficients of half-integral weight cusp forms (for
example, see works by Duke and Iwaniec \cite{Du1, iwaniec3}).

\end{document}